\documentclass[12pt]{article}
\usepackage{amssymb, amsmath, url, graphicx}

\def\3{\subset }
\def\4{\subseteq }
\def\<{\left<}
\def\>{\right>}

\def\bit{\begin{itemize}}
\def\eit{\end{itemize}}
\def\3{\subset }
\def\4{\subseteq }
\def\ov{\overline}

\def\0{\leqno}
\def\a{{\alpha}}
\def\barr{\begin{array}}
\def\earr{\end{array}}
\def\dd{\displaystyle}

\def\Z{{\rlap{$\kern2pt{\rm Z}$}{\rm Z}\,}}

\title{\bf On the factorization numbers\\ of some finite $p$-groups}
\author{Marius T\u arn\u auceanu}
\date{February 17, 2015}

\begin{document}

\maketitle

\begin{abstract}
This note deals with the computation of the factorization number
$F_2(G)$ of a finite group $G$. By using the M\"{o}bius inversion
formula, explicit expressions of $F_2(G)$ are obtained for two
classes of finite abelian groups, improving the results of {\it
Factorization numbers of some finite groups}, Glasgow Math. J.
(2012).
\end{abstract}



\section{Introduction}

Let $G$ be a group, $L(G)$ be the subgroup lattice of $G$ and $H$,
$K$ be two subgroups of $G$. If $G=HK$, then $G$ is said to be
\textit{factorized} by $H$ and $K$ and the expression $G=HK$ is
said to be a \textit{factorization} of $G$. Denote by $F_2(G)$ the
\textit{factorization number} of $G$, that is the number of all
factorizations of $G$.

The starting point for our discussion is given by the paper
\cite{3}, where $F_2(G)$ has been computed for certain classes of
finite groups. The connection between $F_2(G)$ and the subgroup
commutativity degree $sd(G)$ of $G$ (see \cite{5,7}) has been also
established, namely
$$sd(G)=\dd\frac{1}{\mid L(G)\mid^2}\sum_{H\leq G} F_2(H).$$Obviously, by
applying the well-known M\"{o}bius inversion formula to the above
equality, one obtains
$$F_2(G)=\dd\sum_{H\leq G} sd(H)\mid L(H)\mid^2\mu(H,G).\0(1)$$In
particular, if $G$ is abelian, then we have $sd(H)=1$ for all
$H\in L(G)$, and consequently
$$F_2(G)=\dd\sum_{H\leq G} \mid L(H)\mid^2\mu(H,G)=\dd\sum_{H\leq G} \mid
L(G/H)\mid^2\mu(H).\0(2)$$This formula will be used in the
following to calculate the factorization numbers of an elementary
abelian $p$-group and of a rank 2 abelian $p$-group, improving
Theorem 1.2 and Corollary 2.5 of \cite{3}. An interesting
conjecture about the maximum value of $F_2(G)$ on the class of
$p$-groups of the same order will be also presented.
\bigskip

First of all, we recall a theorem due to P. Hall \cite{1} (see
also \cite{2}), that permits us to compute explicitly the
M\"{o}bius function of a finite $p$-group.

\bigskip\noindent{\bf Theorem 1.} {\it Let $G$ be a finite $p$-group of order $p^n$.
Then $\mu(G)=0$ unless $G$ is elementary abelian, in which case we
have $\mu(G)=(-1)^n p^{\,\binom{n}{2}}.$}
\bigskip

In contrast with Theorem 1.2 of \cite{3} that gives only a
recurrence relation satisfied by $F_2(\mathbb{Z}_p^n)$,
$n\in\mathbb{N}$, we are able to determine precise expressions of
these numbers.

\bigskip\noindent{\bf Theorem 2.} {\it We have
$$F_2(\mathbb{Z}_p^n)=\dd\sum_{i=0}^n (-1)^i a_{n,p}(i)\, a_{n-i,p}^2\, p^{\,\binom{i}{2}},\0(3)$$where
$a_{n,p}(i)$ is the number of subgroups of order $p^i$ of
\,$\mathbb{Z}_p^n$, $a_{n,p}$ is the total number of subgroups of
\,$\mathbb{Z}_p^n$, and, by convention, $\binom{i}{2}=0$ for
$i=0,1$.}
\bigskip

Since the numbers $a_{n,p}(i)$, $i=0,1,...,n$, are well-known,
namely
$$a_{n,p}(i)=\frac{(p^n-1)\cdots (p-1)}{(p^i-1)\cdots (p-1)(p^{n-i}-1)\cdots
(p-1)}\,,$$the equality (3) easily leads to the following values
of $F_2(\mathbb{Z}_p^n)$ for $n=1,2,3,4$.\newpage

\bigskip\noindent{\bf Examples.}
\begin{itemize}
\item[{\rm a)}] $F_2(\mathbb{Z}_p)=3$.
\item[{\rm b)}] $F_2(\mathbb{Z}_p^2)=p^2+3p+5$.
\item[{\rm c)}] $F_2(\mathbb{Z}_p^3)=3p^4+4p^3+8p^2+5p+7$.
\item[{\rm d)}] $F_2(\mathbb{Z}_p^4)=p^8+3p^7+9p^6+11p^5+14p^4+15p^3+12p^2+23p+9$.
\end{itemize}
\bigskip

Next we compute the factorization number of a rank 2 abelian
$p$-group.

\bigskip\noindent{\bf Theorem 3.} {\it The factorization number of the finite abelian $p$-group $\mathbb{Z}_{p^{\a_1}}\hspace{-1mm}\times\mathbb{Z}_{p^{\a_2}}$, $\a_1\leq\a_2$, is given by the following equality:
$$F_2(\mathbb{Z}_{p^{\a_1}}\hspace{-1mm}\times\mathbb{Z}_{p^{\a_2}})=\frac{1}{(p-1)^4}\left[(2\a_2-2\a_1+1)p^{2\a_1+4}-(6\a_2-6\a_1+1)p^{2\a_1+3}+\right.$$
$$\left.\hspace{8,5mm}+(6\a_2-6\a_1-1)p^{2\a_1+2}-(2\a_2-2\a_1-1)p^{2\a_1+1}-(2\a_1+2\a_2+3)p^3+\right.$$
$$\left.\hspace{-10mm}+(6\a_1+6\a_2+7)p^2-(6\a_1+6\a_2+5)p+(2\a_1+2\a_2+1)\right].$$}

We remark that Theorem 3 gives a generalization of Corollary 2.5
of \cite{3}. Indeed, by taking $\a_1=1$ and $\a_2=n$ in the above
formula, one obtains:

\bigskip\noindent{\bf Corollary 4.} {\it $F_2(\mathbb{Z}_p\times\mathbb{Z}_{p^n})=(2n-1)p^2+(2n+1)p+(2n+3)\,.$}
\bigskip

Finally, we will focus on the minimum/maximum of $F_2(G)$ when $G$
belongs to the class of $p$-groups of order $p^n$. It is easy to
see that
$$2n+1=F_2(\mathbb{Z}_{p^n})\leq F_2(G).$$For $n\leq 3$ the
greatest value of $F_2(G)$ is obtained for $G\cong\mathbb{Z}_p^n$,
as shows the following result.

\bigskip\noindent{\bf Theorem 5.} {\it Let $G$ be a finite $p$-group of order $p^n$. If $n\leq
3$, then $$F_2(G)\leq F_2(\mathbb{Z}_p^n).$$}

Inspired by Theorem 5, we came up with the following conjecture,
which we also have verified for several $n\geq 4$ and particular
values of $p$.\newpage

\bigskip\noindent{\bf Conjecture 6.} {\it For every finite $p$-group $G$ of order
$p^n$, we have $$F_2(G)\leq F_2(\mathbb{Z}_p^n).$$}

We end our note by indicating a natural problem concerning the
fac\-to\-ri\-za\-tion number of abelian $p$-groups.

\bigskip\noindent{\bf Open problem.} Compute explicitly $F_2(G)$ for
an \textit{arbitrary} finite abelian $p$-group $G$. Given a
positive integer $n$, two partitions $\tau$, $\tau'$ of $n$ and
denoting by $G$, $G'$ the abelian $p$-groups of order $p^n$
induced by $\tau$ and $\tau'$, respectively, is it true that
$F_2(G)\geq F_2(G')$ if and only if $\tau\preceq\,\tau'$ (where
$\preceq$\, denotes the lexicographic order)?

\section{Proofs of the main results}

\bigskip\noindent{\bf Proof of Theorem 2.} By using Theorem 1 in
(2), it follows that
$$F_2(\mathbb{Z}_p^n)=\dd\sum_{H\leq\, \mathbb{Z}_p^n}\mid L(\mathbb{Z}_p^n/H)\mid^2 \mu(H)=\dd\sum_{i=0}^n\dd\sum_{^{\,H\leq\, \mathbb{Z}_p^n}_{\mid H\mid=p^i}}\mid L(\mathbb{Z}_p^n/H)\mid^2 \mu(H)=$$
$$\hspace{19mm}=\dd\sum_{i=0}^n a_{n,p}(i)\mid L(\mathbb{Z}_p^{n-i})\mid^2 (-1)^i\,p^{\,\binom{i}{2}}=\dd\sum_{i=0}^n (-1)^i a_{n,p}(i)\, a_{n-i,p}^2\,
p^{\,\binom{i}{2}},$$as desired.
\hfill\rule{1,5mm}{1,5mm}

\bigskip\noindent{\bf Proof of Theorem 3.} It is well-known that
$G=\mathbb{Z}_{p^{\a_1}}\hspace{-1mm}\times\mathbb{Z}_{p^{\a_2}}$
has a unique elementary abelian subgroup of order $p^2$, say $M$,
and that
$$G/M\cong\mathbb{Z}_{p^{\a_1-1}}\hspace{-1mm}\times\mathbb{Z}_{p^{\a_2-1}}.$$Moreover,
all elementary abelian subgroups of $G$ are contained in
$M$. Denote by $M_i$, $i=1,2,...,p+1$, the minimal subgroups of
$G$. Then every quotient $G/M_i$ is isomorphic to a maximal
subgroup of $G$ and therefore we may assume that
$$G/M_i\cong\mathbb{Z}_{p^{\a_1-1}}\hspace{-1mm}\times\mathbb{Z}_{p^{\a_2}} \mbox{ for } i=1,2,...,p$$and
$$G/M_{p+1}\cong\mathbb{Z}_{p^{\a_1}}\hspace{-1mm}\times\mathbb{Z}_{p^{\a_2-1}}.$$Clearly,
the equality (2) becomes
$$F_2(G)=\hspace{1mm}\mid L(G/M)\mid^2\mu(M)+\dd\sum_{i=1}^{p+1}\mid
L(G/M_i)\mid^2\mu(M_i)+\mid L(G)\mid^2\mu(1),$$in view of Theorem
1. Since by Theorem 2 we have $\mu(M)=\mu(\mathbb{Z}_p^2)=p$,
$\mu(M_i)=\mu(\mathbb{Z}_p)= -1$, for all $i=\ov{1,p+1}$, and
$\mu(1)=1$, one obtains
$$F_2(G)=p\hspace{1mm}\mid L(\mathbb{Z}_{p^{\a_1-1}}\hspace{-1mm}\times\mathbb{Z}_{p^{\a_2-1}})\mid^2-\,p\hspace{1mm}\mid L(\mathbb{Z}_{p^{\a_1-1}}\hspace{-1mm}\times\mathbb{Z}_{p^{\a_2}})\mid^2-\0(4)$$
$$\hspace{3mm}-\hspace{1mm}\mid L(\mathbb{Z}_{p^{\a_1}}\hspace{-1mm}\times\mathbb{Z}_{p^{\a_2-1}})\mid^2+\hspace{1mm}\mid
L(\mathbb{Z}_{p^{\a_1}}\hspace{-1mm}\times\mathbb{Z}_{p^{\a_2}})\mid^2.$$The
total number of subgroups of
$\mathbb{Z}_{p^{\a_1}}\hspace{-1mm}\times\mathbb{Z}_{p^{\a_2}}$
has been computed in Theorem 3.3 of \cite{6}, namely
$$\frac{1}{(p{-}1)^2}\left[(\a_2{-}\a_1{+}1)p^{\a_1{+}2}{-}(\a_2{-}\a_1{-}1)p^{\a_1{+}1}{-}(\a_1{+}\a_2{+}3)p{+}(\a_1{+}\a_2+1)\right].$$Then
the desired formula follows immediately by a direct calculation in
the right side of (4).
\hfill\rule{1,5mm}{1,5mm}

\bigskip\noindent{\bf Proof of Theorem 5.} For $n=2$ we obviously
have
$$F_2(\mathbb{Z}_{p^2})=5<F_2(\mathbb{Z}_p^2)=p^2+3p+5.$$For
$n=3$ it is well-known (see e.g. (4.13), \cite{4}, II) that $G$
can be one of the following groups:
\begin{itemize}
\item[--] $\mathbb{Z}_2^3$, $\mathbb{Z}_2\times\mathbb{Z}_4$, $\mathbb{Z}_8$, $D_8$ and $Q_8$ if $p=2$;
\item[--] $\mathbb{Z}_p^3$, $\mathbb{Z}_p\times\mathbb{Z}_{p^2}$, $\mathbb{Z}_{p^3}$, $M(p^3)=\langle x,y \mid x^{p^2}=y^p=1, y^{-1}xy=x^{p+1}\rangle$ and\newline $E(p^3)=\langle x,y \mid x^p=y^p=[x,y]^p=1, [x,y]\in Z(E(p^3))\rangle$ if $p\geq 3$.
\end{itemize}By using the results in Section 2 of \cite{3}, one
obtains
$$\hspace{-60mm}\mbox{  for }p=2:$$ $$F_2(\mathbb{Z}_2^3)=129>F_2(\mathbb{Z}_2\times\mathbb{Z}_4)=29, F_2(\mathbb{Z}_8)=7, F_2(D_8)=41, F_2(Q_8)=17$$and
$$\hspace{-60mm}\mbox{  for }p\geq 3:$$ $$F_2(\mathbb{Z}_p^3)=3p^4+4p^3+8p^2+5p+7>F_2(\mathbb{Z}_p\times\mathbb{Z}_{p^2})=F_2(M(p^3))=3p^2+5p+7,$$ $$\hspace{-115mm}F_2(\mathbb{Z}_{p^3})=7.$$We
also observe that $E(p^3)$ has $p+1$ elementary abelian subgroups
of order $p^2$, say $M_1, M_2, ..., M_{p+1}$, and that every $M_i$
contains $p+1$ subgroups of order $p$, namely $\Phi(E(p^3))$ and
$M_{ij}$, $j=1,2,...,p$. Then $\hspace{1mm}\mid
L(E(p^3))\mid\,=p^2+2p+4$ and so
$$F_2(E(p^3))<\hspace{1mm}\mid
L(E(p^3))\mid^2\,=p^4+4p^3+12p^2+16p+16.$$On the other hand, we
can easily see that this quantity is less than
$F_2(\mathbb{Z}_p^3)$ for all primes $p\geq 3$, completing the
proof.
\hfill\rule{1,5mm}{1,5mm}

\bigskip\noindent{\bf Remark.} It is clear that an explicit formula for
$F_2(E(p^3))$ cannot be obtained by applying (2), but we are able
to determine it by a direct computation. The factorization pairs
of $E(p^3)$ are:
\begin{itemize}
\item[--] $(1,E(p^3))$, $(E(p^3),1)$;
\item[--] $(M_{ij},M_{i'})\hspace{1mm}\forall\hspace{1mm} i'\neq i$, $(M_{ij},E(p^3))$, $(E(p^3),M_{ij})$, $i=\ov{1,p+1}$, $j=\ov{1,p}$;
\item[--] $(\Phi(E(p^3)),E(p^3))$, $(E(p^3),\Phi(E(p^3)))$;
\item[--] $(M_i,M_{i'j})\hspace{1mm}\forall\hspace{1mm} i'\neq i, j=1,2,...,p$, $(M_i,M_{i'})\hspace{1mm}\forall\hspace{1mm} i'\neq i$, $(M_i,E(p^3))$ and $(M_i,E(p^3))$, $i=\ov{1,p+1}$;
\item[--] $(E(p^3),E(p^3))$.
\end{itemize}Hence $$F_2(E(p^3))=2+p(p+1)(p+2)+2+(p+1)(p^2+p+2)+1=$$
$$\hspace{-38mm}=2p^3+5p^2+5p+7.$$

\bigskip\noindent{\bf Acknowledgements.} The author is grateful to the reviewer for
its remarks which improve the previous version of the paper.

\vspace*{5ex}\small

\hfill
\begin{minipage}[t]{4cm}
Marius T\u arn\u auceanu \\
Faculty of  Mathematics \\
``Al.I. Cuza'' University \\
Ia\c si, Romania \\
e-mail: {\tt tarnauc@uaic.ro}
\end{minipage}

\end{document}